\documentclass[11pt,a4paper]{article}

\usepackage[T1]{fontenc}
\usepackage[utf8]{inputenc}
\usepackage{lmodern}
\usepackage[margin=28mm]{geometry}
\usepackage{microtype}

\usepackage{amsmath,amssymb,amsthm,mathtools}
\IfFileExists{stmaryrd.sty}{\usepackage{stmaryrd}}{%
  \newcommand{\llbracket}{[\![}%
  \newcommand{\rrbracket}{]\!]}%
}
\usepackage{tikz}
\usepackage{hyperref}
\usepackage{bookmark}
\hypersetup{
  colorlinks=true,
  linkcolor=blue,
  filecolor=magenta,
  urlcolor=blue,
  citecolor=blue,
  pdftitle={Cutoff of the interchange process with inhomogeneous conductances},
  pdfsubject={Interchange process, random conductances, and mixing-time cutoff},
  pdfstartview=FitH,
  breaklinks=true
}

\newtheorem{theorema}{Theorem}

\newcommand{\bbE}{\mathbb{E}}
\newcommand{\bbN}{\mathbb{N}}
\newcommand{\bbP}{\mathbb{P}}
\newcommand{\bbR}{\mathbb{R}}
\newcommand{\bbZ}{\mathbb{Z}}
\newcommand{\fg}{\mathbf{g}}
\newcommand{\ga}{\alpha}
\newcommand{\gga}{\gamma}
\newcommand{\gd}{\delta}
\newcommand{\gep}{\varepsilon}
\newcommand{\gD}{\Delta}
\newcommand{\gk}{\kappa}
\newcommand{\gO}{\Omega}
\newcommand{\gl}{\lambda}
\newcommand{\gs}{\sigma}

\newcommand{\cP}{\mathcal{P}}
\newcommand{\cE}{\mathcal{E}}
\newcommand{\cC}{\mathcal{C}}
\newcommand{\cL}{\mathcal{L}}
\newcommand{\cT}{\mathcal{T}}
\newcommand{\cS}{\mathcal{S}}
\newcommand{\bP}{\mathbf{P}}
\newcommand{\bE}{\mathbf{E}}
\newcommand{\ind}{\mathbf{1}}
\newcommand{\fL}{\mathfrak{L}}
\newcommand{\lint}{\llbracket}
\newcommand{\rint}{\rrbracket}
\newcommand{\colonequals}{\coloneqq}
\newcommand{\Var}{\mathrm{Var}}
\newcommand{\TV}{\mathrm{TV}}
\newcommand{\Mix}{\mathrm{mix}}
\newcommand{\tf}{\text{\upshape\scshape f}}
\newcommand{\gap}{\mathrm{gap}}

\DeclareMathOperator{\Proj}{Pr}
\renewcommand{\epsilon}{\varepsilon}
\renewcommand{\tilde}{\widetilde}
\renewcommand{\hat}{\widehat}

\numberwithin{equation}{section}
\theoremstyle{plain}
\newtheorem{theorem}[equation]{Theorem}
\newtheorem{lemma}[equation]{Lemma}
\newtheorem{proposition}[equation]{Proposition}
\newtheorem{corollary}[equation]{Corollary}

\theoremstyle{definition}
\newtheorem{assumption}[equation]{Assumption}
\theoremstyle{remark}
\newtheorem{rem}[equation]{Remark}

\title{Cutoff of the interchange process with inhomogeneous conductances}
\author{Pengfei Tang \thanks{Center for Applied Mathematics and KL-AAGDM, Tianjin University, Tianjin, 300072, China.
		Email: \textsf{pengfei\_tang@tju.edu.cn}.} \and Shangjie Yang\thanks{Department of Applied Mathematics,    Fluminense Federal University, Niterói/RJ, 24210-201, Brasil. 
		Email: \textsf{syang@id.uff.br}.}}
\date{}

\begin{document}

\maketitle

\begin{abstract}
We study the total variation mixing time of the interchange process on the
line segment $\lint 1,N\rint$, where adjacent sites $x$ and $x+1$ exchange
their labels at rate $c^{(N)}(x,x+1)>0$. Under the assumption that the cumulative
resistances $r^{(N)}(a, b)=\sum_{x=a}^{b-1} 1/c^{(N)}(x,x+1)$ are uniformly asymptotic
to the graph distance $b-a$, 
 we show that the chain exhibits  cutoff
at time $(2\pi^2)^{-1}N^2\log N$. In particular, when the conductances are drawn
from an i.i.d.\ (or stationary ergodic) environment with finite mean resistance
$m=\bbE[1/c(1,2)]\in(0,\infty)$, the quenched chains almost surely display cutoff
at time $m(2\pi^2)^{-1}N^2\log N$. 
We adapt Lacoin's strategy to the inhomogeneous conductance setting, by combining
monotonicity and Peres--Winkler censoring with a weighted Dirichlet energy estimate
for the expected height function. By projection, we also deduce cutoff at the same
asymptotic time scale for the simple exclusion process at every macroscopic density
$\rho\in(0,1)$.
\end{abstract}



\section{Introduction}

The interchange process is a fundamental interacting particle system on a finite connected
graph $G=(V,E)$: each vertex initially hosts a particle carrying a distinct label, and
each edge $e\in E$ is equipped with an independent Poisson clock of rate $c(e)>0$ whose
rings transpose the labels at its endpoints. This dynamics naturally couples two other
canonical models on $G$. Declaring $k$ chosen labels to be indistinguishable particles and
the remaining $|V|-k$ labels to be empty sites yields the simple exclusion process with $k$
particles, whereas tracking a single label recovers the continuous-time random walk in the
conductance network $(G,c)$. By Aldous' spectral gap conjecture, established  by Caputo, Liggett, and Richthammer~\cite{caputo2010aldous}, the spectral gaps of
all three Markov chains coincide.

In the present paper, we take $G$ to be the discrete segment $\lint 1,N\rint$, where $\lint x,y\rint\coloneq [x,y ]\cap \mathbb{Z}$, so that
the interchange process becomes the \emph{adjacent transposition shuffle} of a deck of $N$
cards with spatially inhomogeneous (or random) swap rates. In the homogeneous setting
$c\equiv 1$, Wilson~\cite{wilson2004mixing} proved that the total variation mixing time is
sandwiched between $(1+o(1))(2\pi^2)^{-1}N^2\log N$ and $(1+o(1))\pi^{-2}N^2\log N$, and
conjectured that the lower bound is asymptotically sharp. Lacoin~\cite{lacoin2016mixing}
confirmed Wilson's conjecture and established the cutoff phenomenon at time
$(2\pi^2)^{-1}N^2\log N$. Our aim is to prove that this sharp mixing behavior
persists in disordered media---for example, when $(c(x,x+1))_{x\ge 1}$ is a quenched i.i.d.\
sequence with $\bbE[1/c(1,2)]=1$, and more generally whenever the cumulative resistances
satisfy the deterministic law-of-large-numbers condition in Assumption~\ref{assum:resist}
below. In the i.i.d.\ case with unit mean resistance, this condition holds almost surely.

Our proof combines the strategy of~\cite{lacoin2016mixing}---stochastic monotonicity, the
censoring inequality of Peres and Winkler~\cite{peres2013can}, and a decomposition of the
segment into a large but fixed number of blocks---with the spectral asymptotics for the
one-dimensional random conductance model obtained in~\cite{Yang2024LightTail}. A particularly
pleasant feature of Lacoin's argument is that its most delicate ingredient
(Lemma~\ref{lem:skeleton} below) is a  statement concerning the uniform
measure on the symmetric group $\cS_N$ and probability measures whose density is
increasing. Because it involves neither the
dynamics nor the conductances, it carries over verbatim to the inhomogeneous setting.

The disorder enters the proof in only two places: through the asymptotic behavior of the
principal eigenvalue and eigenfunction of the single-particle chain, and through the uniform
exponential decay of the expected height function. The spectral input is furnished directly
by~\cite{Yang2024LightTail}, while for the inhomogeneous discrete heat equation satisfied by
the expected height function (cf.~\cite[Appendix~B]{Yang2024Inhomogeneous}) we derive a
self-contained weighted Dirichlet energy estimate that avoids pointwise bounds on higher
eigenfunctions. Combining these ingredients with a rough coupling bound valid for arbitrary
positive conductances and the inheritance of Assumption~\ref{assum:resist} by macroscopic
subsegments allows us to execute the censoring argument.

In previous work on the inhomogeneous simple exclusion process,
\cite[Theorem~2.6]{Yang2024Inhomogeneous} established cutoff under Assumption~\ref{assum:resist}
together with two additional pointwise hypotheses: a subpolynomial upper bound and a
quantitative lower bound on individual resistances (\cite[Assumptions~2.2 and~2.4]{Yang2024Inhomogeneous}).
By projecting the interchange process onto the exclusion process, our main theorem removes
both pointwise restrictions whenever the particle density $k_N/N$ converges to a limit
$\rho\in(0,1)$, which is  Corollary~\ref{cor:dense}.
\section{Model and results}

\subsection{The interchange process}
For each integer $N\ge 2$, let $\cS_N$ denote the symmetric group on $\lint 1,N\rint$,
equipped with strictly positive edge conductances $(c^{(N)}(x,x+1))_{1\le x<N}$. The interchange
process on $\lint 1,N\rint$ is the continuous-time Markov chain on $\cS_N$ with
generator acting on $f:\cS_N\to\bbR$ by
\begin{equation}\label{generator:interchange}
(\cL_N f)(\gs) \colonequals \sum_{i=1}^{N-1} c^{(N)}(i,i+1)\bigl[f(\gs\circ\tau_{i,i+1})-f(\gs)\bigr],
\end{equation}
where $\tau_{i,j}$ denotes the transposition swapping $i$ and $j$. We interpret $\gs(x)=y$
as saying that site $x$ is occupied by the card labeled $y$. The uniform probability measure
$\mu_N$ on $\cS_N$ is reversible with respect to $\cL_N$ and is its unique
invariant distribution. For $\gs\in\cS_N$, we write $(\eta_t^\gs)_{t\ge 0}$ for
the chain started at $\gs$ and $P_t^\gs$ for the law of $\eta_t^\gs$. 
We define the worst-case distance to equilibrium and the $\gep$-mixing time by
\begin{equation}\label{def:dN}
\begin{aligned}
d_N(t) &\colonequals \max_{\gs\in\cS_N}\|P_t^\gs-\mu_N\|_{\TV},\\
t_{\Mix}^N(\gep) &\colonequals \inf\bigl\{t\ge 0 : d_N(t)\le \gep\bigr\},
\qquad \gep\in(0,1).
\end{aligned}
\end{equation}
where
the total variation distance between two probability measures $\nu_1,\nu_2$ on $\cS_N$
 is $\|\nu_1-\nu_2\|_{\TV}\colonequals \frac{1}{2}\sum_{\gs\in \cS_N}|\nu_1(\gs)-\nu_2(\gs)|$.
We refer to~\cite{LPWMCMT} for background on Markov chain mixing times.
The sequence of chains is said to exhibit \emph{cutoff} if for every $\gep\in(0,1/2)$,
\begin{equation*}
\lim_{N\to\infty}\frac{t^N_{\Mix}(\gep)}{t^N_{\Mix}(1-\gep)} = 1.
\end{equation*}

\begin{figure}[ht]
\centering
\begin{tikzpicture}[x=1cm,y=1cm,font=\small]
  \draw (0,0) -- (1.5,0);
  \draw[densely dotted] (1.5,0) -- (3,0);
  \draw (3,0) -- (5,0);
  \draw[densely dotted] (5,0) -- (6.5,0);
  \draw (6.5,0) -- (8,0);
  \foreach \x/\site/\card in {0/1/{\sigma(1)},1.5/2/{\sigma(2)},3/x/{\sigma(x)},5/{x+1}/{\sigma(x+1)},6.5/{N-1}/{\sigma(N-1)},8/N/{\sigma(N)}} {
    \fill (\x,0) circle (1.6pt);
    \node[below=4pt] at (\x,0) {$\site$};
    \node[above=3pt] at (\x,0) {$\card$};
  }
  \draw[<->] (0,0.8) -- node[above] {$c(1,2)$} (1.5,0.8);
  \draw[<->] (3,0.8) -- node[above] {$c(x,x+1)$} (5,0.8);
  \draw[<->] (6.5,0.8) -- node[above] {$c(N-1,N)$} (8,0.8);
\end{tikzpicture}
\caption{The interchange process on $\lint 1,N\rint$. Each site carries a distinct
label; the labels at adjacent sites $x$ and $x+1$ are exchanged at rate $c(x,x+1)$.}
\label{fig:SEPrandcond}
\end{figure}
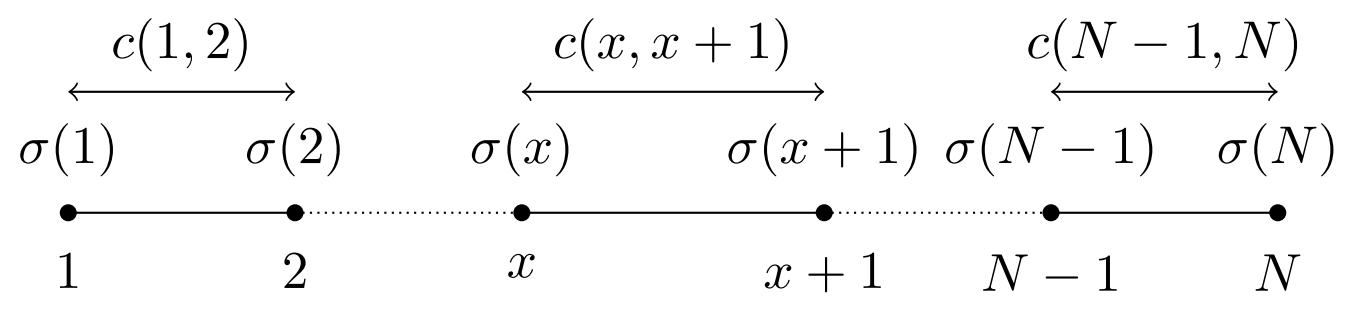

Because the jump rates in~\eqref{generator:interchange} depend only on the edge $\{i,i+1\}$
and not on the current permutation, the trajectory started at $\gs$ can be realized as
$\eta_t^\gs = \gs\circ\pi_t$, where $(\pi_t)_{t\ge 0}$ starts at the identity permutation
$\ind$. Since $\mu_N$ is invariant under left multiplication by $\gs$, it follows that
for every $\gs\in\cS_N$ and $t\ge 0$,
\begin{equation}\label{eq:symmetry}
\|P_t^\gs - \mu_N\|_{\TV} = \|P_t^{\ind} - \mu_N\|_{\TV} = d_N(t).
\end{equation}
Thus, even in an inhomogeneous environment, the total variation distance to equilibrium is
independent of the starting configuration.

\subsection{The simple exclusion process}
For $k\in\lint 1,N-1\rint$, the simple exclusion process with $k$ particles on
$\lint 1,N\rint$ is the Markov chain on
\begin{equation*}
\gO_{N,k} \colonequals \Biggl\{\xi:\lint 1,N\rint \to \{0,1\} \;\Big|\; \sum_{x=1}^N \xi(x) = k\Biggr\}
\end{equation*}
with generator acting on $f:\gO_{N,k}\to\bbR$ by
\begin{equation}\label{generator:SEP}
(\cL_{N,k} f)(\xi) \colonequals \sum_{i=1}^{N-1} c^{(N)}(i,i+1)\bigl[f(\xi\circ\tau_{i,i+1})-f(\xi)\bigr].
\end{equation}
Its unique invariant measure is the uniform distribution $\mu_{N,k}$ on $\gO_{N,k}$, and
we denote its distance to equilibrium and mixing time by $d_{N,k}(t)$ and $t^{N,k}_{\Mix}(\gep)$,
defined analogously to~\eqref{def:dN}. By particle-hole symmetry we may assume $k\le N/2$
whenever convenient.

For each $k\in\lint 1,N-1\rint$, define the projection $\Proj_{N,k}:\cS_N\to\gO_{N,k}$ by
\begin{equation}\label{def:projection}
(\Proj_{N,k}\gs)(x) \colonequals \ind_{\{\gs(x)\le k\}}, \qquad \gs\in\cS_N,\quad x\in\lint 1,N\rint,
\end{equation}
with boundary conventions $\Proj_{N,0}\gs\equiv 0$ and $\Proj_{N,N}\gs\equiv 1$.
In words, the projection $\Proj_{N,k}$ colors the cards labeled $1,\dots,k$ as particles and the remaining
cards as empty sites. Comparing~\eqref{generator:interchange} and~\eqref{generator:SEP}, we see
that $(\Proj_{N,k}\eta_t^\gs)_{t\ge 0}$ is the simple exclusion process started at
$\Proj_{N,k}\gs$, and that $\Proj_{N,k}$ pushes $\mu_N$ forward to $\mu_{N,k}$. Since total
variation distance contracts under projections, we obtain
\begin{equation}\label{TVprojection}
d_{N,k}(t) \le d_N(t), \qquad t\ge 0,\quad k\in\lint 1,N-1\rint.
\end{equation}

\subsection{Spectral gap and Dirichlet form}
The asymptotic rate of exponential convergence to equilibrium is governed by the spectral gap
$\gap_N$, defined as the smallest positive eigenvalue of $-\cL_N$. To express it
variationally, consider the Dirichlet form associated with $\cL_N$: for $f:\cS_N\to\bbR$,
\begin{equation*}
\cE_N(f) \colonequals -\langle f,\cL_N f\rangle_{\mu_N}
= \frac{1}{2}\sum_{\gs\in\cS_N}\sum_{x=1}^{N-1}\mu_N(\gs)c^{(N)}(x,x+1)
\bigl[f(\gs\circ\tau_{x,x+1})-f(\gs)\bigr]^2,
\end{equation*}
where $\langle f,g\rangle_{\mu_N}\colonequals \sum_{\gs\in\cS_N}\mu_N(\gs)f(\gs)g(\gs)$.
Then
\begin{equation}\label{def:gap}
\gap_N = \inf_{\substack{f:\cS_N\to\bbR\\ \Var_{\mu_N}(f)>0}}\frac{\cE_N(f)}{\Var_{\mu_N}(f)},
\end{equation}
where $\Var_{\mu_N}(f)\colonequals \langle f,f\rangle_{\mu_N}-\langle f,\ind\rangle_{\mu_N}^2$; for instance see \cite[Remark~13.8]{LPWMCMT}.
Standard spectral bounds (see~\cite[Chapter~12]{LPWMCMT} in continuous time) yield, for all
$\gep\in(0,1)$,
\begin{equation}\label{relation:gapmixtime}
\frac{1}{\gap_N}\log\frac{1}{2\gep} \le t_{\Mix}^N(\gep) \le \frac{1}{\gap_N}\log\frac{1}{2\gep\mu_{\min}},
\end{equation}
where $\mu_{\min}\colonequals \min_{\gs\in\cS_N}\mu_N(\gs)$, and
\begin{equation}\label{speed:gap}
\lim_{t\to\infty}\frac{1}{t}\log d_N(t) = -\gap_N.
\end{equation}

\subsection{The resistance assumption}
For each edge $\{x,x+1\}$, write $r^{(N)}(x,x+1)\colonequals 1/c^{(N)}(x,x+1)$ for its electrical
resistance, and define the cumulative resistance between $1\le n\le m\le N$ by
$r^{(N)}(n,m)\colonequals \sum_{x=n}^{m-1}r^{(N)}(x,x+1)$ (with the empty sum equal to $0$ when $n=m$).
We suppress the superscript $(N)$ when it is clear from context. Our main results rely on the
following macroscopic regularity condition.

\begin{assumption}\label{assum:resist}
The conductance array satisfies
\begin{equation}\label{LLN}
\gd_N^{(0)} \colonequals \frac{1}{N}\sup_{1\le n<m\le N}\bigl|r^{(N)}(n,m)-(m-n)\bigr| \xrightarrow[N\to\infty]{} 0.
\end{equation}
\end{assumption}

If $(r(x,x+1))_{x\in\bbN}$ is a fixed infinite i.i.d.\ sequence with $\bbE[r(1,2)]=1$,
then~\eqref{LLN} holds almost surely by the strong law of large numbers (see the proof of
Corollary~\ref{cor:iid} below).
Besides Assumption~\ref{assum:resist},  in~\cite[Theorem~2.6]{Yang2024Inhomogeneous} Yang
imposed pointwise subpolynomial  upper and superlogarithmic lower bounds on individual
resistances beyond strict positivity, which are~\cite[Assumptions 2.2, 2.4]{Yang2024Inhomogeneous}.  In our current paper, we only assume  Assumption~\ref{assum:resist} which covers the case:
fix $\ga\in(0,1)$ and set all edge
resistances to $1$ except for one bond of resistance $N^\ga$ and another bond of resistance $(\log N)^{-2}$.  This example does not satisfies either of ~\cite[Assumptions~2.2,~2.4]{Yang2024Inhomogeneous}.

\subsection{Spectral input and main results}
Identifying $\gO_{N,1}$ with $\lint 1,N\rint$, the one-particle generator
$\cL_{N,1}$ acts on $g:\lint 1,N\rint\to\bbR$ by
\[
(\cL_{N,1}g)(x) = c(x,x+1)[g(x+1)-g(x)]\ind_{\{x<N\}} + c(x-1,x)[g(x-1)-g(x)]\ind_{\{x>1\}}.
\]
We first recall its spectral properties from~\cite[Proposition~1.1 and Theorem~1.3]{Yang2024LightTail},
and then relate its nonzero spectrum to that of a weighted discrete Laplacian with Dirichlet boundary conditions (compare~\cite[Appendix~B]{Yang2024Inhomogeneous}).

\begin{theorema}\label{th:A}
For any choice of positive conductances, $-\cL_{N,1}$ has simple eigenvalues
$0=\gl_0<\gl_1<\cdots<\gl_{N-1}$. Let $g_i$ be the eigenfunction satisfying
$\cL_{N,1}g_i = -\gl_i g_i$ normalized by $g_i(1)=1$. Then $g_1$ is strictly
decreasing. Furthermore, under Assumption~\ref{assum:resist}, we have
\begin{equation}\label{eigenvalue:asym}
\lim_{N\to\infty}\frac{N^2\gl_1}{\pi^2} = 1,
\end{equation}
and, with $h_1(x)\colonequals \cos(\pi(x-1/2)/N)$,
\begin{equation}\label{eigenfunction:asym}
\lim_{N\to\infty}\max_{x\in\lint 1,N\rint}|g_1(x)-h_1(x)| = 0, 	\qquad
\lim_{N \to \infty} \sup_{x \in \lint 2, N \rint}\left| N(c \nabla g_1)(x)-N(\nabla h_1)(x)\right|\;=\;0\,,
\end{equation}
where $(c \nabla f)(x)\colonequals c(x-1,x)\left[f(x)-f(x-1)\right]$. 
\end{theorema}
Only the principal eigenvalue and eigenfunction asymptotics in~\eqref{eigenvalue:asym}
and~\eqref{eigenfunction:asym} will be needed in the sequel.

\begin{theorema}\label{th:B}
Set $S_N\colonequals \sum_{x=1}^{N-1}r(x,x+1)$ and for
$x\in\lint 1,N-1\rint$, $\varpi_N(x)\colonequals r(x,x+1)/S_N$. On the space of functions $u:\lint 0,N\rint\to\bbR$
satisfying $u(0)=u(N)=0$, define the operator
\begin{equation}\label{def:heightoperator}
(A_N u)(x) \colonequals c(x,x+1)\bigl[2u(x)-u(x+1)-u(x-1)\bigr], \qquad x\in\lint 1,N-1\rint.
\end{equation}
Then $A_N$ is self-adjoint and positive definite on $L^2(\lint 1,N-1\rint,\varpi_N)$,
and admits an orthonormal eigenbasis $(\fg_i)_{1\le i\le N-1}$ with eigenvalues
$0<\gk_1<\cdots<\gk_{N-1}$:
\begin{equation}\label{def:heateigen}
A_N\fg_i = \gk_i\fg_i, \qquad \fg_i(0)=\fg_i(N)=0.
\end{equation}
Moreover, for arbitrary positive conductances, we have
\begin{equation}\label{eq:spectra}
\gk_i = \gl_i, \qquad 1\le i\le N-1.
\end{equation}
In particular, under Assumption~\ref{assum:resist}, we have $\lim_{N\to\infty} N^2\gk_1 = \pi^2$.
\end{theorema}

\begin{proof}
Since $\varpi_N(x)c(x,x+1)=1/S_N$, summation by parts gives, for any functions $u,v$ vanishing at
$0$ and $N$,
\begin{equation}\label{heat:dirichlet}
\begin{aligned}
\langle u,A_N v\rangle_{\varpi_N}
&= \frac{1}{S_N}\sum_{x=1}^{N-1}u(x)\bigl[2v(x)-v(x+1)-v(x-1)\bigr]\\
&= \frac{1}{S_N}\sum_{x=0}^{N-1}\bigl[u(x+1)-u(x)\bigr]\bigl[v(x+1)-v(x)\bigr],
\end{aligned}
\end{equation}
where the second equality follows by splitting the bracket as $[v(x)-v(x-1)]-[v(x+1)-v(x)]$,
shifting the index in the first sum, and using $u(0)=u(N)=0$. This expression is symmetric in
$u$ and $v$, and strictly positive when $u=v\not\equiv 0$. The finite-dimensional spectral theorem (cf.\ \cite[Theorem~7.29]{Axler2024})
therefore provides the orthonormal eigenbasis $(\fg_i)_{1\le i\le N-1}$ with positive eigenvalues.

To identify the spectrum, set $H_i(x)\colonequals \sum_{z=1}^x g_i(z)$ for $1\le i\le N-1$ and
$x\in\lint 0,N\rint$ (with $H_i(0)=0$). Since $\cL_{N,1}$ is reversible with respect to the uniform measure, the eigenfunction $g_i$, corresponding to a nonzero eigenvalue,  is orthogonal to the constant
eigenfunction $g_0\equiv 1$. Hence we  have $H_i(N)=0$.
Summing $(\cL_{N,1}g_i)(z) = -\gl_i g_i(z)$ over $z\in\lint 1,x\rint$ yields
\[
c(x,x+1)\bigl[g_i(x+1)-g_i(x)\bigr] = -\gl_i H_i(x), \qquad x\in\lint 1,N-1\rint.
\]
Since $g_i(x+1)-g_i(x) = [H_i(x+1)-H_i(x)]-[H_i(x)-H_i(x-1)]$, we obtain $A_N H_i = \gl_i H_i$.
Because $H_i(1)=g_i(1)=1$, each $H_i$ is a nontrivial eigenfunction of $A_N$ with eigenvalue $\gl_i$.
Since $0<\gl_1<\cdots<\gl_{N-1}$ are $N-1$ distinct eigenvalues on an $(N-1)$-dimensional
space, they exhaust the spectrum of $A_N$, proving~\eqref{eq:spectra}. The final assertion then
follows from~\eqref{eigenvalue:asym}.
\end{proof}

We are now ready to state our main result.

\begin{theorem}\label{thm:main}
Under Assumption~\ref{assum:resist}, for every $\gep\in(0,1)$, we have
\begin{equation}\label{eq:main}
\lim_{N\to\infty}\frac{2\pi^2 t_{\Mix}^N(\gep)}{N^2\log N} = 1.
\end{equation}
In particular, the sequence of interchange processes displays total variation cutoff at time
$(2\pi^2)^{-1}N^2\log N$.
\end{theorem}

\begin{corollary}\label{cor:iid}
Let $(c(x,x+1))_{x\in\bbN}$ be a fixed infinite sequence of i.i.d.\ strictly positive
conductances under a probability law $\bbP$, and assume that $m\colonequals \bbE[1/c(1,2)]\in(0,\infty)$.
For each $N\ge 2$, endow $\lint 1,N\rint$ with the first $N-1$ conductances. Then, for
$\bbP$-almost every environment and every $\gep\in(0,1)$, we have
\begin{equation}\label{eq:iid}
\lim_{N\to\infty}\frac{2\pi^2 t^N_{\Mix}(\gep)}{m N^2\log N} = 1,
\end{equation}
so the quenched chains exhibit cutoff at time $\frac{m}{2\pi^2}N^2\log N$.
\end{corollary}

\begin{proof}
Suppose first that $m=1$, and set $R_j\colonequals \sum_{x=1}^j r(x,x+1)$ for $j\ge 1$, with $R_0\colonequals 0$.
By the strong law of large numbers, $(R_j-j)/j\to 0$ almost surely as $j\to\infty$. Fix any $\eta>0$;
almost surely there exists $j_0$ such that $|R_j-j|\le \eta j\le \eta N$ for all $j_0<j<N$, while
$N^{-1}\max_{0\le j\le j_0}|R_j-j|\to 0$ as $N\to\infty$. Since $\eta>0$ is arbitrary, it follows
that $N^{-1}\max_{0\le j<N}|R_j-j|\to 0$ almost surely. Writing
\[
r(n,m')-(m'-n) = \bigl(R_{m'-1}-(m'-1)\bigr)-\bigl(R_{n-1}-(n-1)\bigr)
\]
shows that Assumption~\ref{assum:resist} holds almost surely, and Theorem~\ref{thm:main}
yields~\eqref{eq:iid} with $m=1$. For general $m\in(0,\infty)$, replacing $c$ by $mc$ rescales the
resistances to have unit mean and accelerates the continuous-time dynamics by a factor of $m$, which
divides the mixing time by $m$.
\end{proof}

\begin{rem}[Stationary ergodic environments]\label{rem:ergodic}
Corollary~\ref{cor:iid} extends without change to any stationary and ergodic sequence of strictly
positive conductances $(c(x,x+1))_{x\in\bbN}$ with $m=\bbE[1/c(1,2)]\in(0,\infty)$.
Indeed, Birkhoff's ergodic theorem implies $R_n/n\to m$ almost surely, and the same argument in the proof of Corollary~\ref{cor:iid}  verifies Assumption~\ref{assum:resist} after dividing the resistances by $m$.
\end{rem}

\begin{corollary}\label{cor:dense}
Under Assumption~\ref{assum:resist}, suppose that $k_N/N\to\rho\in(0,1)$ as $N\to\infty$. Then for
every $\gep\in(0,1)$, we have 
\begin{equation}\label{eq:dense}
\lim_{N\to\infty}\frac{2\pi^2 t^{N,k_N}_{\Mix}(\gep)}{N^2\log N} = 1.
\end{equation}
The same limit holds almost surely in the i.i.d.\ setting of Corollary~\ref{cor:iid} and the
stationary ergodic setting of Remark~\ref{rem:ergodic}, with $N^2\log N$ replaced by $mN^2\log N$.
\end{corollary}

\begin{proof}
The asymptotic upper bound follows immediately from~\eqref{TVprojection} and Theorem~\ref{thm:main},
while the matching lower bound is given by Proposition~\ref{lem:denselower} together with~\eqref{eigenvalue:asym}.
The i.i.d.\ and stationary ergodic statements follow by the same  rescaling as in Corollary~\ref{cor:iid}
and Remark~\ref{rem:ergodic}.
\end{proof}

\medskip\noindent\textbf{Organization of the paper.}
Section~\ref{sec:prelim} gathers the preliminaries needed for the analysis: the height partial order on
$\cS_N$, the graphical construction, the equality $\gap_N=\gap_{N,1}$, the censoring and
FKG inequalities, a rough coupling bound valid for arbitrary positive conductances, and the decay
estimate for the expected height function. Section~\ref{sec:lowerbound} presents the lower bound on the mixing time for the simple exclusion process and deduces
the lower bound for the interchange process. Section~\ref{sec:upperbound} establishes the matching
upper bound, completing the proof of Theorem~\ref{thm:main}.

\section{Preliminaries}\label{sec:prelim}

\subsection{Partial orders and height functions}\label{subsec:order}
Following~\cite[Section~3.1]{lacoin2016mixing}, we encode each permutation $\gs\in\cS_N$
by a centered height surface $\tilde{\gs}:\lint 0,N\rint^2\to\bbR$ defined by
\begin{equation}\label{def:surface}
\tilde{\gs}(x,y) \colonequals \sum_{z=1}^x \ind_{\{\gs(z)\le y\}} - \frac{xy}{N}, \qquad (x,y)\in\lint 0,N\rint^2,
\end{equation}
where the linear term $xy/N$ ensures that $\tilde{\gs}(x,y)$ has zero mean under $\mu_N$.
This encoding is injective because for all $(x,y)\in\lint 1,N\rint^2$,
\begin{equation*}
\tilde{\gs}(x,y)-\tilde{\gs}(x,y-1)-\tilde{\gs}(x-1,y)+\tilde{\gs}(x-1,y-1)+\frac{1}{N} = \ind_{\{\gs(x)=y\}}.
\end{equation*}
Identifying $\gs$ with $\tilde{\gs}$, we endow $\cS_N$ with the coordinatewise partial order
\begin{equation}\label{def:order}
\gs \preceq \gs' \quad \Longleftrightarrow \quad \tilde{\gs}(x,y) \le \tilde{\gs}'(x,y) \quad \text{for all } (x,y)\in\lint 0,N\rint^2.
\end{equation}
Since $\tilde{\gs}(x,y)\le \min(x,y)-xy/N = \tilde{\ind}(x,y)$, the identity permutation
$\ind$ is the unique maximal element of $(\cS_N,\preceq)$, while the reversal
$\gs_{\min}(x)\colonequals N+1-x$ is the unique minimal element.

Similarly, for $\xi\in\gO_{N,k}$, we define its height function $h^\xi:\lint 0,N\rint\to\bbR$
by $h^\xi(0)\colonequals 0$ and
\begin{equation}\label{def:height}
h^\xi(x) \colonequals \sum_{y=1}^x \xi(y) - \frac{k}{N}x, \qquad x\in\lint 1,N\rint,
\end{equation}
and equip $\gO_{N,k}$ with the order $\xi\preceq\xi'\iff h^\xi(x)\le h^{\xi'}(x)$ for all
$x\in\lint 0,N\rint$. Its maximal and minimal configurations are $\wedge_{N,k}\colonequals \ind_{\{1\le x\le k\}}$
and $\vee_{N,k}\colonequals \ind_{\{N-k<x\le N\}}$, respectively. Comparing~\eqref{def:projection}
and~\eqref{def:surface}, we observe that
\begin{equation}\label{eq:surfaceheight}
h^{\Proj_{N,k}\gs}(x) = \tilde{\gs}(x,k), \qquad x\in\lint 0,N\rint,
\end{equation}
so the $k$-th horizontal slice of $\tilde{\gs}$ is precisely the height function of the projected
exclusion configuration, and $\Proj_{N,k}$ is order-preserving. Moreover, every $\xi\in\gO_{N,k}$
satisfies the elementary uniform bound
\begin{equation}\label{bound:height}
\max_{x\in\lint 0,N\rint}|h^\xi(x)| \le \frac{k(N-k)}{N} \le \min(k,N-k).
\end{equation}

\subsection{A graphical construction}\label{subsec:graphconstr}
We construct a grand coupling of $(\eta_t^\gs)_{t\ge 0}$ across all initial permutations
$\gs\in\cS_N$ as in~\cite[Section~3.2]{lacoin2016mixing}, adapting only the Poisson clock
intensities to the inhomogeneous conductances. For each $x\in\lint 1,N-1\rint$, let
$\cT^x=(\cT_n^x)_{n\ge 1}$ be a Poisson point process on $\bbR_+$ of intensity
$2c(x,x+1)$, and let $(U_n^x)_{x,n}$ be an array of i.i.d.\ $\mathrm{Bernoulli}(1/2)$ random variables,
all mutually independent. Denote their joint law and expectation by $\bP$ and $\bE$ respectively.

For each $\gs\in\cS_N$, the c\`adl\`ag trajectory $(\eta_t^\gs)_{t\ge 0}$ starts at
$\eta_0^\gs=\gs$, remains constant between clock rings, and at each arrival time $t=\cT_n^x$
updates $\eta_{t^-}^\gs$ as follows:
\begin{itemize}
\item if $U_n^x=1$ and $\eta_{t^-}^\gs(x)>\eta_{t^-}^\gs(x+1)$, swap the labels at $x$ and $x+1$ (\emph{sort});
\item if $U_n^x=0$ and $\eta_{t^-}^\gs(x)<\eta_{t^-}^\gs(x+1)$, swap the labels at $x$ and $x+1$ (\emph{reverse sort});
\item otherwise, leave $\eta_{t^-}^\gs$ unchanged.
\end{itemize}
From any permutation, exactly one of the two values of $U_n^x$ triggers the transposition $\tau_{x,x+1}$,
so each adjacent transposition occurs at rate $2c(x,x+1)\cdot\frac{1}{2}=c(x,x+1)$, matching the generator~\eqref{generator:interchange}.
We refer to the action at time $\cT_n^x$ as an \emph{update at coordinate $x$}.

In terms of the surface $\tilde{\gs}$, an update at coordinate $x$ acts independently on each
level $y$: if $\tilde{\eta}_{t^-}^\gs(\cdot,y)$ has a local minimum at $x$ and $U_n^x=1$, then
$\tilde{\eta}_t^\gs(x,y)=\tilde{\eta}_{t^-}^\gs(x,y)+1$; symmetrically, if it has a local
maximum at $x$ and $U_n^x=0$, then $\tilde{\eta}_t^\gs(x,y)=\tilde{\eta}_{t^-}^\gs(x,y)-1$;
in all other cases $\tilde{\eta}^\gs(\cdot,y)$ is unchanged. Because this local update rule does
not depend on the clock intensities, the proof of~\cite[Proposition~3.1]{lacoin2016mixing} applies
verbatim and yields the following monotonicity property.

\begin{proposition}\label{prop:attractive}
For all $\gs,\gs'\in\cS_N$,
\begin{equation*}
\gs \preceq \gs' \quad \Longrightarrow \quad \bP\bigl[\eta_t^\gs \preceq \eta_t^{\gs'} \text{ for all } t\ge 0\bigr] = 1.
\end{equation*}
Similarly, projecting by $\Proj_{N,k}$ yields an order-preserving grand coupling of the simple
exclusion trajectories $(\eta_t^\xi)_{t\ge 0}$, $\xi\in\gO_{N,k}$.
\end{proposition}

\subsection{Spectral gap invariance across particle numbers}
The following statement combines Aldous' spectral gap theorem~\cite[Theorem~1.1]{caputo2010aldous}
with the strict monotonicity of the principal one-particle eigenfunction from Theorem~\ref{th:A};
see also~\cite[Proposition~3.2]{Yang2024Inhomogeneous}.

\begin{proposition}\label{prop:gapSEP}
For arbitrary positive conductances,
\begin{equation}\label{gapinterchange}
\gap_N = \gap_{N,k} = \gap_{N,1} = \gl_1, \qquad k\in\lint 1,N-1\rint.
\end{equation}
Moreover, the linear functional $\tf(\xi)\colonequals \sum_{x=1}^N \xi(x)g_1(x)$ on $\gO_{N,k}$
satisfies $\cL_{N,k}\tf = -\gl_1\tf$. Setting
$\gd_{\min}\colonequals \min_{1\le x<N}[g_1(x)-g_1(x+1)]>0$, the function $\tf$ is increasing for the partial order $\preceq$ on $\gO_{N,k}$
and satisfies $\tf(\xi')-\tf(\xi)\ge \gd_{\min}$ whenever $\xi\preceq\xi'$ and $\xi\ne\xi'$.
\end{proposition}

\begin{proof}
All the  statements are proved in~\cite[Proposition~3.2]{Yang2024Inhomogeneous},  except that $\gap_N = \gl_1$, which  follows from~\cite[Theorem~1.1]{caputo2010aldous}. Indeed the proof in~\cite[Proposition~3.2]{Yang2024Inhomogeneous} together with the projections $(\Proj_{N,k})_k$ implies~\eqref{gapinterchange}.

\end{proof}

In particular, under Assumption~\ref{assum:resist}, Theorem~\ref{th:A} implies $\lim_{N\to\infty} N^2\gap_N = \pi^2$.

\subsection{Monotone measures, censoring, and correlation inequalities}\label{subsec:censoring}
A function $f:\cS_N\to\bbR$ is called \emph{increasing} if $\gs\preceq\gs'$ implies
$f(\gs)\le f(\gs')$. A probability measure $\nu$ on $\cS_N$ is called \emph{increasing}
if $\gs\mapsto\nu(\gs)$ (equivalently, the density $\nu/\mu_N$) is increasing. Recall that
$\nu$ stochastically dominates $\nu'$, denoted by $\nu\succeq\nu'$, if $\nu(f)\ge \nu'(f)$ for every increasing
function $f$. Here  $\nu(f)$ denotes the expectation of $f$ under the measure $\nu$, i.e., $\nu(f)\coloneq \sum_{\sigma\in\cS_N}f(\sigma)\nu(\sigma)$. By the FKG inequality (Proposition~\ref{prop:FKG} below) applied to $\nu/\mu_N$ and $f$,
every increasing probability measure $\nu$ satisfies $$\nu(f)=\mu_N((\nu/\mu_N)f)\ge \mu_N(\nu/\mu_N)\mu_N(f)=\mu_N(f),$$
and therefore stochastically dominates $\mu_N$.

A \emph{censoring scheme} is a c\`adl\`ag map $\cC:\bbR_+\to \cP(\lint 1,N-1\rint)$,
where $\cP(\lint 1,N-1\rint)$ is the power set of $\lint 1,N-1\rint$. The censored
dynamics associated with $\cC$ is obtained from the graphical construction of
Subsection~\ref{subsec:graphconstr} by executing the update at coordinate $x$ at time $t=\cT_n^x$
if and only if $x\notin\cC(t)$, and suppressing it otherwise. We denote by $P_t^{\nu,\cC}$
the marginal law at time $t$ of the censored dynamics started from an initial distribution $\nu$.

\begin{proposition}\label{prop:censoring}
Let $\nu$ be an increasing probability measure on $\cS_N$ and let $\cC$ be a censoring
scheme. Then $P_t^\nu$ and $P_t^{\nu,\cC}$ are increasing for every $t\ge 0$, and
\begin{equation}\label{censoringineq}
\|P_t^\nu - \mu_N\|_{\TV} \le \|P_t^{\nu,\cC} - \mu_N\|_{\TV}.
\end{equation}
\end{proposition}

\begin{proof}
This is~\cite[Propositions~A.1 and~3.6]{lacoin2016mixing}, which adapts the censoring inequality of
Peres and Winkler~\cite[Theorem~1.1]{peres2013can} to $\cS_N$. The argument carries over
to arbitrary positive conductances without change. Indeed, for $x\in\lint 1,N-1\rint$ and
$\gs\in\cS_N$, let $\gs_x^\bullet\colonequals \{\gs,\gs\circ\tau_{x,x+1}\}$.
A single update at coordinate $x$ maps a probability measure $\nu$ to $\theta_x\nu$, where
\begin{equation*}
(\theta_x\nu)(\gs) \colonequals \frac{1}{2}\nu(\gs_x^\bullet) = \frac{1}{2}\bigl[\nu(\gs)+\nu(\gs\circ\tau_{x,x+1})\bigr]\,.
\end{equation*}
Thus the update  operator $\theta_x$ does not depend on the conductances. The proof in~\cite[Appendix~A.2]{lacoin2016mixing}
relies solely on two combinatorial properties of $\theta_x$ established in~\cite[Lemma~A.2]{lacoin2016mixing}---namely,
that $\theta_x$ preserves the class of increasing measures and that $\nu\succeq\theta_x(\nu)$ whenever
$\nu$ is increasing---followed by averaging over the random sequence of update coordinates. Since the
averaging step holds for any distribution of update times, the conclusion remains valid in the
inhomogeneous setting.
\end{proof}

We next recall the block projections and static correlation inequalities from~\cite{lacoin2016mixing}.
Given an integer parameter $K\in\lint 1,N\rint$, define the cut points and block lengths
\begin{equation}\label{def:xi}
x_i \colonequals \lceil iN/K\rceil, \quad 0\le i\le K; \qquad \gD x_i \colonequals x_i - x_{i-1}, \quad 1\le i\le K.
\end{equation}
We partition both positions and labels into $K$ intervals:
\begin{equation}\label{def:blocks}
I_i \colonequals \lint x_{i-1}+1,x_i\rint \quad\text{(positions)}, \qquad
J_i \colonequals \lint x_{i-1}+1,x_i\rint \quad\text{(labels)}, \qquad 1\le i\le K.
\end{equation}
Although $I_i$ and $J_i$ are the same subsets of $\lint 1,N\rint $, we use distinct notation to distinguish blocks of positions from blocks of labels. 
The \emph{semi-skeleton} $\hat{\gs}$ and \emph{skeleton} $\bar{\gs}$ of $\gs\in\cS_N$
are the restrictions of $\tilde{\gs}$ given by
\begin{equation}\label{def:skeleton}
\hat{\gs}(x,j) \colonequals \tilde{\gs}(x,x_j), \qquad
\bar{\gs}(i,j) \colonequals \tilde{\gs}(x_i,x_j), \qquad (x,i,j)\in\lint 0,N\rint\times\lint 0,K\rint^2,
\end{equation}
with image sets $\hat{\cS}_N$ and $\bar{\cS}_N$ equipped with the coordinatewise order.
The semi-skeleton records which label block $J_j$ occupies each position $x$, since
\[
\hat{\gs}(x,j) - \hat{\gs}(x-1,j) + \frac{x_j}{N} = \ind_{\{\gs(x)\le x_j\}}\,.
\]
As a useful consequence, $\hat\sigma=\hat{\ind}$ if and only if $\sigma(x)$ and $x$ belong to the same label block for every $x\in\lint 1,N\rint$, or equivalently, $\sigma(J_i)=J_i$
 for every $i\in\lint 1,K\rint$. In contrast,  the skeleton retains only the total count of labels from each $J_j$ lying in each position
block $I_i$ (see~\eqref{eq:blockcounts}). For a probability measure $\nu$ on $\cS_N$, we denote
by $\hat{\nu}$, $\bar{\nu}$, and $\bar{\nu}_{i,j}$ its push-forwards under $\gs\mapsto\hat{\gs}$,
$\gs\mapsto\bar{\gs}$, and $\gs\mapsto\bar{\gs}(i,j)$, respectively.

\begin{proposition}[FKG inequality on $\cS_N$, {\cite[Proposition~3.5]{lacoin2016mixing}}]\label{prop:FKG}
For any pair of increasing functions $f,g:\cS_N\to\bbR$, we have 
\[
\mu_N\bigl(fg\bigr) \ge \mu_N(f)\mu_N(g).
\]
\end{proposition}

\begin{proposition}[Preservation of monotonicity under projection, {\cite[Proposition~3.8]{lacoin2016mixing}}]\label{prop:projection}
If $\nu$ is an increasing probability measure on $\cS_N$, then $\bar{\nu}/\bar{\mu}_N$ is
increasing on $\bar{\cS}_N$, and for each $(i,j)\in\lint 1,K\rint^2$ the ratio
$\bar{\nu}_{i,j}/\bar{\mu}_{N,i,j}$ is increasing on the support of $\bar{\mu}_{N,i,j}$.
\end{proposition}

Finally, let $\tilde{\cS}_N\colonequals \{\gs\in\cS_N : \hat{\gs}=\hat{\ind}\}$
be the subgroup of permutations leaving every label block $J_i$ invariant, which is isomorphic to
$\prod_{i=1}^K\cS_{\gD x_i}$. For any probability measure $\nu$ on $\cS_N$, define
its left-averaged measure
\begin{equation}\label{def:tildenu}
\tilde{\nu}(\gs) \colonequals \frac{1}{|\tilde{\cS}_N|}\sum_{\pi\in\tilde{\cS}_N}\nu(\pi\circ\gs)
= \frac{1}{|\tilde{\cS}_N|}\hat{\nu}(\hat{\gs}),
\end{equation}
where the second equality holds because the fibers of $\gs\mapsto\hat{\gs}$ are precisely the
right cosets $\tilde{\cS}_N\gs$.

Averaging on the left commutes with the censored evolution: writing $A\nu\colonequals \tilde{\nu}$,
associativity of group multiplication gives, for every coordinate $x$,
\[
(A\theta_x\nu)(\gs)
= \frac{1}{2|\tilde{\cS}_N|}\sum_{\pi\in\tilde{\cS}_N}\bigl[\nu(\pi\circ\gs)+\nu(\pi\circ\gs\circ\tau_{x,x+1})\bigr]
= (\theta_x A\nu)(\gs).
\]
Conditioning on the finite sequence of uncensored clock rings in $[0,t]$, the censored dynamics is a
finite composition of the operators $\theta_x$; commuting $A$ through this composition and taking
expectation over the clocks yields
\begin{equation}\label{eq:commute}
\widetilde{P_t^{\nu,\cC}} = P_t^{\tilde{\nu},\cC}, \qquad t\ge 0, \quad \text{for every censoring scheme }\cC.
\end{equation}

\begin{lemma}[{\cite[Lemma~4.3]{lacoin2016mixing}}]\label{lem:decomposition}
For every probability measure $\nu$ on $\cS_N$, we have $\|\tilde{\nu}-\mu_N\|_{\TV}=\|\hat{\nu}-\hat{\mu}_N\|_{\TV}$,
and consequently
\begin{equation}\label{decomposition}
\|\nu-\mu_N\|_{\TV} \le \|\hat{\nu}-\hat{\mu}_N\|_{\TV} + \|\nu-\tilde{\nu}\|_{\TV}.
\end{equation}
\end{lemma}

\subsection{A rough upper bound on the distance to equilibrium}\label{subsec:rough}
We next establish an exponential bound on the distance to equilibrium,  valid for arbitrary positive conductances and
analogous to~\cite[Proposition~4.2]{lacoin2016mixing} (see also~\cite[Proposition~C.2]{Yang2024Inhomogeneous}).
For a finite interval $J\subset\bbZ$ of cardinality $\ell\ge 2$ with positive edge conductances,
write $d_J(t)$ for the distance to equilibrium of the interchange process on $J$ and
$0=\gl_0(J)<\gl_1(J)<\cdots<\gl_{\ell-1}(J)$ for the eigenvalues of $-\cL_{J,1}$, the one-particle chain on $J$.

\begin{proposition}\label{prop:rough}
For every finite segment $J$ of length $\ell\ge 2$, every assignment of positive conductances on $J$,
and every $t\ge 0$, we have
\begin{equation}\label{roughbound}
d_J(t) \le  \ell^2 e^{-\gl_1(J)t}.
\end{equation}
\end{proposition}

\begin{proof}
By translation we may assume $J=\lint 1,\ell\rint$. Since $\mu_\ell$ is invariant, the
triangle inequality gives
\begin{equation}\label{rough:step1}
d_J(t) \le \max_{\gs,\gs'\in\cS_\ell}\|P_t^\gs - P_t^{\gs'}\|_{\TV}.
\end{equation}
Fix $\gs,\gs'\in\cS_\ell$ and construct a Markovian coupling $(\eta_t^\gs,\eta_t^{\gs'})_{t\ge 0}$
as in~\cite[Proposition~6.5]{lacoin2016mixing} and~\cite[Appendix~C]{Yang2024Inhomogeneous}: for each
$x\in\lint 1,\ell-1\rint$,
\begin{itemize}
\item if $\eta_t^\gs(x)\ne \eta_t^{\gs'}(x)$ and $\eta_t^\gs(x+1)\ne \eta_t^{\gs'}(x+1)$,
then the transition $\eta\mapsto \eta\circ \tau_{x,x+1}$ occurs independently for the two processes, with rate $c(x,x+1)$ each;
\item otherwise (when $\eta_t^\gs$ and $\eta_t^{\gs'}$ agree at $x$ or at $x+1$), the transition $\eta\mapsto \eta\circ \tau_{x,x+1}$  is performed simultaneously by the two processes, with rate  $c(x,x+1)$.
\end{itemize}
Each marginal is an interchange process with generator~\eqref{generator:interchange}. For each label
$i\in\lint 1,\ell\rint$, let $X_t^i\colonequals (\eta_t^\gs)^{-1}(i)$ and
$Y_t^i\colonequals (\eta_t^{\gs'})^{-1}(i)$ denote its positions in the two configurations.
Whenever $X_t^i=Y_t^i=z$, i.e.\ both chains agree at $z$, so any move along an incident edge moves $X^i$
and $Y^i$ together. Before they meet, each coordinate jumps along its incident edges at the
one-particle rates: if $|X_t^i-Y_t^i|\ge 2$, the sets of edges incident to $X_t^i$ and $Y_t^i$ are
disjoint and therefore driven by independent clocks; if $\{X_t^i,Y_t^i\}=\{x,x+1\}$, both endpoints
of $\{x,x+1\}$ are mismatched, so the two coordinates cross $\{x,x+1\}$ independently. Hence
$(X_t^i,Y_t^i)_{t\ge 0}$ evolves as two independent one-particle random walks until their first
meeting time and remains on the diagonal since then.

Let $\fL$ denote the generator of this coupled two-particle walk and $\bP_{x,y}$ its
law started at $(x,y)$. By the union bound (see~\cite[Lemma~2.1]{lacoin2016mixing}),
\begin{equation}\label{rough:step2}
\|P_t^\gs - P_t^{\gs'}\|_{\TV}
\le \bP\bigl[\exists\,i\in\lint 1,\ell\rint : X_t^i\ne Y_t^i\bigr]
\le \ell\max_{x\ne y}\bP_{x,y}[X_t\ne Y_t].
\end{equation}
To estimate the meeting probability, set
\[
D \colonequals \{(x,y)\in\lint 1,\ell\rint^2 : x<y\}, \qquad
T \colonequals \inf\{s\ge 0 : X_s = Y_s\},
\]
and consider the killed transition kernel $G_t^*(a,b)\colonequals \bP_a[(X_t,Y_t)=b,\ t<T]$
for $a,b\in D$. Since nearest-neighbor walks on $\lint 1,\ell\rint$ cannot cross without
meeting, a process started in $D$ remains in $D$ up to time $T$. Moreover, the killed generator is
symmetric with respect to counting measure on $D$ because every allowed jump and its reverse carry
the same edge conductance.

Let $(\tilde{g}_i)_{0\le i\le \ell-1}$ be an eigenbasis of $-\cL_{\ell,1}$ orthonormal in
$L^2(\lint 1,\ell\rint,\mu_{\ell,1})$, so that
\begin{equation}\label{rough:orthonormal}
\sum_{x=1}^\ell \tilde{g}_i(x)\tilde{g}_j(x) = \ell\,\gd_{i,j}.
\end{equation}
For $0\le i<j\le \ell-1$, define the antisymmetric functions
\begin{equation}\label{def:uij}
u_{i,j}(x,y) \colonequals \tilde{g}_i(x)\tilde{g}_j(y) - \tilde{g}_i(y)\tilde{g}_j(x).
\end{equation}
Each $u_{i,j}$ vanishes on the diagonal and satisfies $\fL u_{i,j}=-(\gl_i+\gl_j)u_{i,j}$
off the diagonal (see \cite[Lemma~C.3]{Yang2024Inhomogeneous}). Using antisymmetry and~\eqref{rough:orthonormal}, we compute for $i<j$ and $p<q$:
\begin{equation}\label{rough:norm}
\sum_{1\le x<y\le \ell}u_{i,j}(x,y)u_{p,q}(x,y)
= \frac{1}{2}\sum_{x,y=1}^\ell u_{i,j}(x,y)u_{p,q}(x,y)
= \ell^2(\gd_{i,p}\gd_{j,q}-\gd_{i,q}\gd_{j,p})
= \ell^2\gd_{(i,j),(p,q)}.
\end{equation}
Since there are $\binom{\ell}{2}=|D|$ pairs $(i,j)$ with $0\le i<j\le \ell-1$, the normalized family
$(u_{i,j}/\ell)_{i<j}$ forms an orthonormal eigenbasis of $\ell^2(D)$. Consequently,
\begin{equation}\label{rough:kernel}
G_t^*\bigl((x_0,y_0),(x,y)\bigr) = \sum_{0\le i<j\le \ell-1}\frac{u_{i,j}(x_0,y_0)u_{i,j}(x,y)}{\ell^2}e^{-(\gl_i+\gl_j)t}.
\end{equation}
Because $\gl_i+\gl_j\ge \gl_0+\gl_1=\gl_1$ for all $i<j$, we obtain
$\|G_t^*\|_{\ell^2(D)\to\ell^2(D)} = e^{-\gl_1 t}$. Writing $\gd_a$ for the Dirac mass at
$a\in D$, the Cauchy--Schwarz inequality gives
\begin{equation}\label{rough:step3}
\bP_a(T>t) = \langle \gd_a, G_t^*\ind\rangle_{\ell^2(D)}
\le \|\gd_a\|_2 \|G_t^*\|_{2\to 2} \|\ind\|_2
= \sqrt{\binom{\ell}{2}}\,e^{-\gl_1 t}.
\end{equation}
The same estimate holds for $x_0>y_0$ by symmetry. Combining~\eqref{rough:step1}, \eqref{rough:step2},
and~\eqref{rough:step3} yields~\eqref{roughbound}.
\end{proof}

\subsection{Inheritance of the resistance condition by blocks}
In Section~\ref{sec:upperbound} we partition $\lint 1,N\rint$ into a fixed number $K$ of blocks
$I_1,\dots,I_K$. The next lemma ensures that each block inherits Assumption~\ref{assum:resist}.

\begin{lemma}\label{lem:blocks}
Fix $K\in\bbN$ and let $(I_i)_{1\le i\le K}$ be the position blocks from~\eqref{def:blocks}.
If $(c(x,x+1))_{1\le x<N}$ satisfies Assumption~\ref{assum:resist}, then for each $i\in\lint 1,K\rint$
the restricted conductances $(c(x,x+1))_{x_{i-1}<x<x_i}$ satisfy Assumption~\ref{assum:resist} with
$N$ replaced by $\gD x_i$. In particular, by Theorem~\ref{th:A},
\begin{equation}\label{gapblocks}
\lim_{N\to\infty}(\gD x_i)^2\gl_1(I_i) = \pi^2, \qquad i\in\lint 1,K\rint.
\end{equation}
\end{lemma}

\begin{proof}
Set $\ell\colonequals \gD x_i$. By~\eqref{def:xi}, $N/K-1\le \ell\le N/K+1$, so $N\le (K+1)\ell$
whenever $N\ge K(K+1)$. For any $x_{i-1}<n<m\le x_i$, we have $|r(n,m)-(m-n)|\le N\gd_N^{(0)}$, whence
\begin{equation*}
\frac{1}{\ell}\sup_{x_{i-1}<n<m\le x_i}\bigl|r(n,m)-(m-n)\bigr|
\le \frac{N}{\ell}\gd_N^{(0)} \le (K+1)\gd_N^{(0)} \xrightarrow[N\to\infty]{} 0.
\end{equation*}
Because $K$ is fixed, the convergence is uniform over $i\in\lint 1,K\rint$, and~\eqref{gapblocks}
follows from~\eqref{eigenvalue:asym}.
\end{proof}

\subsection{Decay of the expected height function}\label{subsec:heat}
By~\eqref{def:surface} and~\eqref{eq:surfaceheight}, for each fixed $y\in\lint 1,N-1\rint$
the process $(\tilde{\eta}_t^\gs(\cdot,y))_{t\ge 0}$ is the height function of the $y$-particle
simple exclusion process. As noted in~\cite[Appendix~B]{Yang2024Inhomogeneous}, its expectation
$f(t,x)\colonequals \bE[\tilde{\eta}_t^\gs(x,y)]$ satisfies the inhomogeneous discrete heat
equation
\begin{equation}\label{heatequation}
\begin{cases}
\partial_t f(t,x) = c(x,x+1)\bigl[f(t,x+1)+f(t,x-1)-2f(t,x)\bigr], & x\in\lint 1,N-1\rint,\\
f(t,0) = f(t,N) = 0, & t\ge 0,\\
f(0,x) = \tilde{\gs}(x,y).
\end{cases}
\end{equation}
Indeed, at coordinate $x$, $\tilde{\eta}_t(\cdot,y)$ jumps upward by $1$ at rate $c(x,x+1)$ when it
has a local minimum at $x$, downward by $1$ at rate $c(x,x+1)$ when it has a local maximum at $x$,
and does not jump otherwise, so that the drift is $c(x,x+1)[\tilde{\eta}_t(x+1,y)+\tilde{\eta}_t(x-1,y)-2\tilde{\eta}_t(x,y)]$.
In terms of the spectral decomposition of $A_N$ from Theorem~\ref{th:B}, the solution of~\eqref{heatequation} is
\begin{equation}\label{heat:solution}
\bE\bigl[\tilde{\eta}_t^\gs(x,y)\bigr]
= \sum_{i=1}^{N-1}e^{-\gk_i t}\fg_i(x)\bigl\langle \fg_i,\tilde{\gs}(\cdot,y)\bigr\rangle_{\varpi_N}.
\end{equation}
The following lemma extends~\cite[Lemma~4.1]{lacoin2016mixing} (see also~\cite[Lemma~C.1]{Yang2024Inhomogeneous})
without requiring pointwise bounds on the eigenfunctions $\fg_i$.

\begin{lemma}\label{lem:heat}
	Recall $S_N\colonequals \sum_{x=1}^{N-1}r(x,x+1)$ from Theorem~\ref{th:B}.
For arbitrary positive conductances, every $\gs\in\cS_N$, every $y\in\lint 0,N\rint$,
and every $t\ge (2\gk_1)^{-1}$, we have
\begin{equation}\label{heat:general}
\max_{x\in\lint 0,N\rint}\bigl|\bE\bigl[\tilde{\eta}_t^\gs(x,y)\bigr]\bigr|
\le \sqrt{N S_N \gk_1}\,\min(y,N-y)e^{-\gk_1 t}.
\end{equation}
Furthermore, 
under Assumption~\ref{assum:resist}, for all sufficiently large $N$ (with a threshold independent of
$\gs$, $y$, and $t$), we have
\begin{equation}\label{bound:heat}
\max_{x\in\lint 0,N\rint}\bigl|\bE\bigl[\tilde{\eta}_t^\gs(x,y)\bigr]\bigr|
\le 2\pi\min(y,N-y)e^{-\gk_1 t}, \qquad t\ge (2\gk_1)^{-1}.
\end{equation}
\end{lemma}

\begin{proof}
The cases $y\in\{0,N\}$ are trivial since $\tilde{\gs}(\cdot,y)\equiv 0$. Fix $y\in\lint 1,N-1\rint$
and write $f_t(x)\colonequals \bE[\tilde{\eta}_t^\gs(x,y)]$. By \eqref{heatequation} and the definition of $A_N$ in~\eqref{def:heightoperator}, we have $\partial_t f_t=-A_Nf_t$, hence  $f_t = e^{-tA_N}f_0$ with
$f_0 = \tilde{\gs}(\cdot,y)$. Expanding $f_0 = \sum_{i=1}^{N-1}a_i\fg_i$ in the orthonormal
basis of Theorem~\ref{th:B}, Parseval's identity gives $\sum_{i=1}^{N-1}a_i^2 = \|f_0\|_{L^2(\varpi_N)}^2$.
Since $\frac{d}{dz}(ze^{-2tz})=e^{-2tz}(1-2tz)\le 0$ for $z\ge \gk_1$ and $t\ge (2\gk_1)^{-1}$,
we have $\gk_i e^{-2\gk_i t}\le \gk_1 e^{-2\gk_1 t}$ for every $i\ge 1$, and hence
\begin{equation}\label{heat:energydecay}
\langle f_t,A_N f_t\rangle_{\varpi_N}
= \sum_{i=1}^{N-1}\gk_i a_i^2 e^{-2\gk_i t}
\le \gk_1 e^{-2\gk_1 t}\|f_0\|_{L^2(\varpi_N)}^2.
\end{equation}
On the other hand, any function $u$ vanishing at $0$ and $N$ satisfies
$u(x)=\sum_{j=0}^{x-1}[u(j+1)-u(j)]$ for $x\in\lint 1,N-1\rint$. Applying Cauchy--Schwarz
and~\eqref{heat:dirichlet} yields the bound
\begin{equation}\label{heat:sobolev}
\|u\|_\infty^2 \le N\sum_{j=0}^{N-1}\bigl[u(j+1)-u(j)\bigr]^2 = N S_N \langle u,A_N u\rangle_{\varpi_N}.
\end{equation}
Note that $f_t(0)=f_t(N)=0$. Applying \eqref{heat:sobolev} to $u=f_t$ and then using \eqref{heat:energydecay}, we obtain
\[
\|f_t\|_\infty^2
\leq NS_N\langle f_t,A_Nf_t\rangle_{\varpi_N}
\leq NS_N\gk_1e^{-2\gk_1 t}\|f_0\|_{L^2(\varpi_N)}^2.
\]
Since $\sum_{x=1}^{N-1}\varpi_N(x)=1$, \eqref{bound:height} gives 
\[
\|f_0\|_{L^2(\varpi_N)}\le \|f_0\|_\infty \le \min(y,N-y)\,.
\]
Substituting this bound and taking square roots establishes~\eqref{heat:general}.
Finally, under Assumption~\ref{assum:resist}, we have $S_N/N = r(1,N)/N\to 1$ while $N^2\gk_1\to \pi^2$
by Theorem~\ref{th:B}. Therefore $N S_N \gk_1 \to \pi^2$, so $\sqrt{N S_N \gk_1}\le 2\pi$ for all
sufficiently large $N$, proving~\eqref{bound:heat}.
\end{proof}

\section{Lower bound on the mixing time}\label{sec:lowerbound}

This section is a corollary of~\eqref{TVprojection} with $k_N=\lfloor N/2\rfloor$ and \cite[Proposition 4.1]{Yang2024Inhomogeneous}, which does not require \cite[Assumptions~2.2,~2.4]{Yang2024Inhomogeneous}. We explain the  mechanism for the positive particle density case behind \cite[Proposition 4.1]{Yang2024Inhomogeneous}: one uses the
eigenfunction $\tf(\xi)=\sum_x \xi(x)g_1(x)$ of $\cL_{N,k}$ (see Proposition \ref{prop:gapSEP}) as a
distinguishing statistic. The process $e^{\gl_1 (t-t_0)}\tf(\eta_t)$ is a martingale, so that
$\bE\left[\tf(\eta^{\wedge}_{t_0})\right]=e^{-\gl_1t_0}\tf(\wedge)$ is of order $e^{-\gl_1 t_0}k$, while
the martingale bracket, controlled by means of \eqref{eigenfunction:asym}, shows that
the fluctuations of $\tf$ under both $P^{\wedge}_{t_0}$ and $\mu_{N,k}$ are $O(\sqrt k)$. The two are
of the same order precisely at $t_0 = \frac{1}{2\gl_1}\log k$.

\begin{proposition}\label{prop:lower}
Under Assumption~\ref{assum:resist}, for every $\gep\in(0,1)$, we have 
\begin{equation}\label{lowerbound}
\liminf_{N\to\infty}\frac{2\pi^2 t^N_{\Mix}(\gep)}{N^2\log N} \ge 1.
\end{equation}
\end{proposition}

\begin{proposition}{(\cite[Proposition 4.1]{Yang2024Inhomogeneous})}\label{lem:denselower}	
Suppose Assumption~\ref{assum:resist} holds and there exists $\gga>0$ such that  $  q_N:=\min(N-k_N, k_N) \ge  (\log N)^{1+\gga}$.  For all $\gep \in (0, 1)$ and for all $N$ sufficiently large we have
\begin{equation}\label{lower:dense}
	t_{\Mix}^{N,k}(\gep) \;\ge\; \frac{1}{2\gl_1} \log q_N-\frac{ C(\gep)}{\gl_1}\,.
\end{equation}	
\end{proposition}

\begin{proof}[Proof of Proposition~\ref{prop:lower}]
Applying~\eqref{TVprojection} and Proposition~\ref{lem:denselower} with $k_N=\lfloor N/2\rfloor$ gives
\[
t^N_{\Mix}(\gep) \ge t^{N,k_N}_{\Mix}(\gep) \ge \frac{\log N}{2\gl_1} - \frac{C(\gep)}{\gl_1}.
\]
Dividing by $N^2\log N/(2\pi^2)$ and invoking~\eqref{eigenvalue:asym} yields~\eqref{lowerbound}.
\end{proof}

\section{Upper bound on the mixing time}\label{sec:upperbound}

We now turn to the upper bound, which forms the core of the paper. Together with
Proposition~\ref{prop:lower} and Theorem~\ref{th:B}, the following proposition immediately
implies Theorem~\ref{thm:main}.

\begin{proposition}\label{prop:upper}
Under Assumption~\ref{assum:resist}, for every $\gep\in(0,1)$, $\gd\in(0,1/300)$, and
for all sufficiently large $N$ we have
\begin{equation}\label{upperbound}
d_N\left((1+\gd)\frac{\log N}{2\gk_1}\right) \le \gep.
\end{equation}
Consequently, $\limsup_{N\to\infty} 2\pi^2 t^N_{\Mix}(\gep)/(N^2\log N) \le 1$.
\end{proposition}

\subsection{The censoring scheme and proof strategy}
Fix $\gep\in(0,1)$ and $\gd\in(0,1/300)$ throughout this section, and set
\begin{equation}\label{def:K}
K \colonequals \lceil 1/\gd\rceil,
\end{equation}
so that $K\ge 300$ and $\gd K\ge 1$. All thresholds in $N$ below depend only on the fixed
parameters $(\gep,\gd)$ and the conductance sequence, and are chosen uniformly over the
$K$ blocks and $K^2$ skeleton coordinates.

Recall the cut points $x_i=\lceil iN/K\rceil$ from~\eqref{def:xi} and the blocks $I_i,J_i$
from~\eqref{def:blocks}. Define the three times
\begin{equation}\label{def:times}
t_1 \colonequals \frac{\gd}{3}\cdot\frac{\log N}{2\gk_1}, \qquad
t_2 \colonequals \left(1+\frac{2\gd}{3}\right)\frac{\log N}{2\gk_1}, \qquad
t_3 \colonequals (1+\gd)\frac{\log N}{2\gk_1},
\end{equation}
and consider the censoring scheme
\begin{equation}\label{def:censoring}
\cC(t) \colonequals
\begin{cases}
\{x_1,\dots,x_{K-1}\} & \text{if } t\in [0,t_1)\cup[t_2,t_3),\\
\emptyset & \text{if } t\in [t_1,t_2)\cup[t_3,\infty).
\end{cases}
\end{equation}
During $[0,t_1)$ and $[t_2,t_3)$, exchanges across the block boundaries $x_1,\dots,x_{K-1}$ are
forbidden, so conditioned on the configuration at the start of the interval, the dynamics evolves
independently within each of the $K$ blocks $I_1,\dots,I_K$; during $[t_1,t_2)$ and $[t_3,\infty)$,
all edges evolve freely.

Write $\nu_t\colonequals P_t^{\ind,\cC}$ for the marginal law at time $t$ of this
censored chain started from the identity permutation. Since $\ind$ is the maximal element of
$(\cS_N,\preceq)$, the Dirac mass $\gd_{\ind}$ is an increasing probability measure.
Combining~\eqref{eq:symmetry} with Proposition~\ref{prop:censoring} therefore gives
\begin{equation}\label{upper:censoring}
\nu_t \text{ is increasing for every } t\ge 0, \qquad\text{and}\qquad
d_N(t_3) = \|P_{t_3}^{\ind}-\mu_N\|_{\TV} \le \|\nu_{t_3}-\mu_N\|_{\TV}.
\end{equation}
Thus it suffices to show that at time $t_3$, the total variation distance between  equilibrium and the censored chain is at most $\gep$.
Following~\cite[Section~5]{lacoin2016mixing}, we break the evolution of $\nu_t$ into three stages:
by time $t_1$, the relative order of cards within each block is close to equilibrium (Proposition~\ref{prop:step1});
by time $t_2$, the skeleton is close to equilibrium (Proposition~\ref{prop:step2}); and during
$[t_2,t_3]$,  mixing within the blocks brings the semi-skeleton close to equilibrium (Proposition~\ref{prop:step3}).

\begin{proposition}\label{prop:step1}
For all sufficiently large $N$, $\|\nu_t-\tilde{\nu}_t\|_{\TV}\le \gep/3$ for every $t\ge t_1$.
\end{proposition}

\begin{proposition}\label{prop:step2}
For all sufficiently large $N$, $\|\bar{\nu}_{t_2}-\bar{\mu}_N\|_{\TV}\le \gep/3$.
\end{proposition}

\begin{proposition}\label{prop:step3}
For all sufficiently large $N$, $\|\hat{\nu}_{t_3}-\hat{\mu}_N\|_{\TV}\le 2\gep/3$.
\end{proposition}

\begin{proof}[Proof of Proposition~\ref{prop:upper}]
Combining~\eqref{upper:censoring}, Lemma~\ref{lem:decomposition}, and Propositions~\ref{prop:step1}
and~\ref{prop:step3}, we obtain for all sufficiently large $N$
\begin{equation*}
d_N(t_3) \le \|\nu_{t_3}-\mu_N\|_{\TV} \le \|\hat{\nu}_{t_3}-\hat{\mu}_N\|_{\TV} + \|\nu_{t_3}-\tilde{\nu}_{t_3}\|_{\TV}
\le \frac{2\gep}{3} + \frac{\gep}{3} = \gep. \qedhere
\end{equation*}
\end{proof}

Before proving Propositions~\ref{prop:step1}--\ref{prop:step3}, we establish a bound on mixing within each block that will be used during both censored intervals $[0,t_1)$ and $[t_2,t_3)$, each of 
duration $t_1$.

\begin{lemma}\label{lem:blockmixing}
For all sufficiently large $N$,
\begin{equation}\label{blockmixing}
\sum_{i=1}^K d_{I_i}(t_1) \le \frac{\gep}{3}.
\end{equation}
\end{lemma}

\begin{proof}
Set $\ell_i\colonequals \gD x_i$. For $N\ge 2K$, every block has length $2\le \ell_i\le N/K+1\le 2N/K$.
Since $K$ is fixed, Lemma~\ref{lem:blocks} yields $\min_{1\le i\le K}\ell_i^2\gl_1(I_i)\ge \pi^2/2$
for all sufficiently large $N$, so that
\[
\gl_1(I_i) \ge \frac{\pi^2}{2\ell_i^2} \ge \frac{\pi^2 K^2}{8N^2}, \qquad i\in\lint 1,K\rint.
\]
Furthermore, Theorem~\ref{th:B} ensures that $\pi^2/(N^2\gk_1)\ge 24/25$ for all large $N$. Using
$\gd K\ge 1$ and $K\ge 300$, we obtain
\[
t_1\gl_1(I_i)
\ge \frac{\gd K^2}{48}\frac{\pi^2}{N^2\gk_1}\log N
\ge \frac{\gd K^2}{50}\log N
\ge \frac{K}{50}\log N
\ge 6\log N.
\]
Applying Proposition~\ref{prop:rough} on each block $I_i$ gives
\begin{equation*}
\sum_{i=1}^K d_{I_i}(t_1)
\le \sum_{i=1}^K \ell_i^2 e^{-\gl_1(I_i)t_1}
\le K\left(\frac{2N}{K}\right)^2 N^{-6}
= \frac{4}{K N^4} \xrightarrow[N\to\infty]{} 0,
\end{equation*}
which proves~\eqref{blockmixing}.
\end{proof}

\subsection{Proof of Proposition~\ref{prop:step1}}\label{subsec:step1}
\begin{proof}
Because updates across $\{x_1,\dots,x_{K-1}\}$ are censored during $[0,t_1)$ and the chain starts
at $\ind$, no card leaves its initial block up to time $t_1$: for every $t\le t_1$,
$\eta_t^{\ind,\cC}$ maps each $I_i$ onto $I_i=J_i$, meaning
$\eta_t^{\ind,\cC}\in\tilde{\cS}_N$ almost surely. Hence $\hat{\nu}_t = \gd_{\hat{\ind}}$,
and by~\eqref{def:tildenu},
\begin{equation}\label{step1:tildenu}
\tilde{\nu}_t = \tilde{\gd}_{\ind} = \text{the uniform probability measure on }\tilde{\cS}_N, \qquad t\le t_1.
\end{equation}
On $[0,t_1)$, the censored chain decomposes as a product of $K$ independent interchange processes
living on $I_1,\dots,I_K$, each started at the identity. Writing $\nu_t^i$ for the law of the $i$-th
block and $\mu^i$ for the uniform measure on $\cS_{I_i}$, subadditivity of total variation
distance for product measures and Lemma~\ref{lem:blockmixing} yield
\begin{equation}\label{step1:bound}
\|\nu_{t_1}-\tilde{\gd}_{\ind}\|_{\TV}
\le \sum_{i=1}^K \|\nu_{t_1}^i-\mu^i\|_{\TV}
\le \sum_{i=1}^K d_{I_i}(t_1) \le \frac{\gep}{3}
\end{equation}
for all sufficiently large $N$, establishing the bound at $t=t_1$.

Finally, by~\eqref{eq:commute}, $\tilde{\nu}_t = P_t^{\tilde{\gd}_{\ind},\cC}$ is
the time-$t$ marginal of the very same censored dynamics started from $\tilde{\gd}_{\ind}$.
For any $t\ge t_1$, both $\nu_t$ and $\tilde{\nu}_t$ are obtained by propagating $\nu_{t_1}$ and
$\tilde{\nu}_{t_1}$ through the same Markov transition kernel from $t_1$ to $t$, so
$\|\nu_t-\tilde{\nu}_t\|_{\TV}\le \|\nu_{t_1}-\tilde{\nu}_{t_1}\|_{\TV}\le \gep/3$ by the
contraction property of Markov kernels (for instance see \cite[Exercise~4.2]{LPWMCMT}).
\end{proof}

\subsection{Proof of Proposition~\ref{prop:step2}}\label{subsec:step2}
To control the skeleton at time $t_2$, define the cumulative skeleton volume
\begin{equation}\label{def:volume}
v(\bar{\gs}) \colonequals \sum_{i,j=1}^{K-1}\bar{\gs}(i,j).
\end{equation}
The key  input from~\cite{lacoin2016mixing} states that if a measure with increasing density
has expected skeleton volume much smaller than $\sqrt{N}$ (the scale of equilibrium fluctuations),
then its skeleton is close to $\bar{\mu}_N$ in total variation.

\begin{lemma}[{\cite[Lemma~5.5]{lacoin2016mixing}}]\label{lem:skeleton}
Fix $K\ge 2$. For every $\gep>0$, there exists $\eta=\eta(K,\gep)>0$ such that for
all sufficiently large $N$ and every probability measure $\nu$ on $\cS_N$ with increasing
density $\nu/\mu_N$,
\begin{equation}\label{skeletonimplication}
\nu\bigl(v(\bar{\gs})\bigr) \le \eta\sqrt{N} \quad \Longrightarrow \quad \|\bar{\nu}-\bar{\mu}_N\|_{\TV} \le \gep,
\end{equation}
where $\nu(v(\bar{\gs}))$ denotes the expectation of $v(\bar{\gs})$ under $\nu$.
\end{lemma}

As emphasized in the introduction, Lemma~\ref{lem:skeleton} involves only the uniform measure
$\mu_N$ and the partial order~\eqref{def:order}; since neither the dynamics nor the conductances
enter its statement, \cite[Lemma~5.5]{lacoin2016mixing} applies verbatim.

\begin{proof}[Proof of Proposition~\ref{prop:step2}]
During $[t_1,t_2)$, no edges are censored, so conditioned on $\eta_{t_1}^{\ind,\cC}=\gs$,
the process $(\eta_{t_1+s}^{\ind,\cC})_{0\le s\le t_2-t_1}$ evolves as the unconstrained
interchange process started at $\gs$. By~\eqref{def:times},
\[
t_2 - t_1 = \left(1+\frac{\gd}{3}\right)\frac{\log N}{2\gk_1} \ge \frac{1}{2\gk_1}
\]
for all $N\ge 3$. Applying Lemma~\ref{lem:heat} to the process started from $\sigma$ and then averaging over $\sigma$ with law  $\nu_{t_1}$, we obtain
 for every $(i,j)\in\lint 1,K-1\rint^2$ and all sufficiently large $N$,
\[
\bigl|\nu_{t_2}\bigl(\bar{\gs}(i,j)\bigr)\bigr|
= \Bigl|\bE\bigl[\tilde{\eta}_{t_2}^{\ind,\cC}(x_i,x_j)\bigr]\Bigr|
\le 2\pi\min(x_j,N-x_j)e^{-\gk_1(t_2-t_1)}
\le 2\pi N \cdot N^{-1/2-\gd/6}.
\]
Summing over $i,j\in\lint 1,K-1\rint$ yields
\begin{equation}\label{step2:volume}
\nu_{t_2}\bigl(v(\bar{\gs})\bigr) \le 2\pi(K-1)^2 N^{-\gd/6}\sqrt{N}.
\end{equation}
Since $K$ and $\gd$ are fixed, we have $2\pi(K-1)^2 N^{-\gd/6}\le \eta(K,\gep/3)$
for all sufficiently large $N$. Because $\nu_{t_2}$ is increasing by~\eqref{upper:censoring},
Lemma~\ref{lem:skeleton} implies $\|\bar{\nu}_{t_2}-\bar{\mu}_N\|_{\TV}\le \gep/3$.
\end{proof}

\subsection{Proof of Proposition~\ref{prop:step3}}\label{subsec:step3}
\begin{proof}
	Write $\gs_t\colonequals \eta_t^{\ind,\cC}$ and follow~\cite[Section~5.5]{lacoin2016mixing}.
	The process runs without censoring on $[t_1,t_2)$, so irreducibility gives
	$\nu_{t_2}(\rho)>0$ for every $\rho\in\cS_N$.
	During $[t_2,t_3)$, censoring blocks exchanges across $\{x_1,\dots,x_{K-1}\}$, and since almost
	surely no Poisson clock rings at the deterministic time $t_3$, each image set $\gs_t(I_i)$
	remains constant on $[t_2,t_3]$. Conditioned on $\gs_{t_2}$, the dynamics inside the $K$ blocks
	$I_1,\dots,I_K$ evolve as independent interchange processes over a time interval of length
	$t_3-t_2=t_1$, with stationary distribution
	\begin{equation}\label{eq:mu_rho}
	\mu_{\gs_{t_2}} \colonequals \mu_N\bigl(\,\cdot\mid \gs(I_i)=\gs_{t_2}(I_i) \text{ for all } i\in\lint 1,K\rint\bigr).
	\end{equation}
	By subadditivity of total variation distance for product measures and Lemma~\ref{lem:blockmixing}, for
	every realization of $\gs_{t_2}$ and all sufficiently large $N$, we have 
	\begin{equation}\label{step3:blocks}
	\bigl\|\bP[\gs_{t_3}\in\cdot\mid \gs_{t_2}] - \mu_{\gs_{t_2}}\bigr\|_{\TV}
	\le \sum_{i=1}^K d_{I_i}(t_1) \le \frac{\gep}{3}.
	\end{equation}
	
	For a fixed permutation $\rho\in\cS_N$, write $\mu_\rho$ for the measure
	in~\eqref{eq:mu_rho} when $\gs_{t_2}=\rho$.
	We claim that the distribution $\hat{\mu}_\rho$ of the semi-skeleton under $\mu_\rho$ satisfies
	\begin{equation}\label{eq:conditionalequilibrium}
	\hat{\mu}_\rho=\hat{\mu}_N(\cdot\mid\bar{\gs}=\bar{\rho}).
	\end{equation}
	
	For $i,j\in\lint 1,K\rint$, let $n_{ij}(\gs)\colonequals|\gs(I_i)\cap J_j|$.
	By~\eqref{def:surface} and~\eqref{def:skeleton}, taking differences at
	successive cut points gives
	\begin{equation}\label{eq:blockcounts}
	n_{ij}(\gs)=\bar{\gs}(i,j)-\bar{\gs}(i-1,j)-\bar{\gs}(i,j-1)
	+\bar{\gs}(i-1,j-1)+\frac{\gD x_i\gD x_j}{N}.
	\end{equation}
	Conversely, summing the block counts up to each pair of cut points and
	subtracting the centering term recovers the skeleton.
	
	Associate to each permutation $\gs$ the $K$-tuple of words $(w_1,\dots,w_K)$,
	where $w_i\in\lint 1,K\rint^{I_i}$ and $w_i(x)=j$ if and only if
	$\gs(x)\in J_j$. The resulting tuples are in bijective correspondence
	with the semi-skeletons $\hat{\gs}$.
	Call a word $w_i$ $\rho$-admissible if it contains exactly
	$n_{ij}(\rho)$ occurrences of $j$ for every $j$. Call a tuple
	$(w_1,\dots,w_K)$ $\rho$-admissible if each $w_i$ is $\rho$-admissible. Under $\mu_\rho$ the label sets
	$\rho(I_i)$ fix these counts, and under
	$\mu_N(\cdot\mid\bar{\gs}=\bar{\rho})$ the skeleton fixes them.
	
	By the bijective correspondence above, to prove~\eqref{eq:conditionalequilibrium}
	it suffices to show that $(w_1,\dots,w_K)$ is uniform on the set of
	$\rho$-admissible tuples under both $\mu_\rho$ and
	$\mu_N(\cdot\mid\bar{\gs}=\bar{\rho})$.
	
	Under $\mu_\rho$, the restrictions $\gs|_{I_i}$ are independent uniform
	bijections $I_i\to\rho(I_i)$. To realize a fixed $\rho$-admissible word $w_i$,
	assign the labels in $\rho(I_i)\cap J_j$ bijectively to the positions
	in $I_i$ marked $j$. This gives $n_{ij}(\rho)!$ choices for each $j$,
	so exactly $\prod_{j=1}^K n_{ij}(\rho)!$ bijections induce $w_i$.
	The count is independent of the word, so the tuple is uniform on the
	set of $\rho$-admissible tuples.
	
	Under $\mu_N(\cdot\mid\bar{\gs}=\bar{\rho})$, fix a $\rho$-admissible tuple.
	Across all the words, there are $\sum_i n_{ij}(\rho)=|J_j|=\gD x_j$
	positions marked $j$. Assigning the labels in $J_j$ bijectively to these
	positions gives $(\gD x_j)!$ choices for each $j$. Every resulting
	permutation has skeleton $\bar\rho$, since the block counts determine
	the skeleton. Thus exactly $\prod_{j=1}^K(\gD x_j)!$ permutations
	induce the fixed tuple. This count is also independent of the tuple,
	so under $\mu_N(\cdot\mid\bar{\gs}=\bar{\rho})$ the tuple $(w_1,\dots,w_K)$
	is also uniform on the set of $\rho$-admissible tuples.
	This proves~\eqref{eq:conditionalequilibrium}.
	
	Projecting the two laws in~\eqref{step3:blocks} to the semi-skeleton gives
	\[
	\bigl\|\bP[\hat{\gs}_{t_3}\in\cdot\mid\gs_{t_2}=\rho]
	-\hat{\mu}_\rho\bigr\|_{\TV}\le\frac{\gep}{3}
	\]
	for every $\rho\in\cS_N$. Fix $\xi\in\bar{\cS}_N$. For each $\rho$ with
	$\bar{\rho}=\xi$,
	substituting~\eqref{eq:conditionalequilibrium} into the preceding inequality gives
	\[
	\bigl\|\bP[\hat{\gs}_{t_3}\in\cdot\mid\gs_{t_2}=\rho]
	-\hat{\mu}_N(\cdot\mid\bar{\gs}=\xi)\bigr\|_{\TV}
	\le\frac{\gep}{3}.
	\]
	The events $\{\bar{\gs}_{t_2}=\xi\}$ and $\{\bar{\gs}_{t_3}=\xi\}$ coincide
	almost surely by preservation of the skeleton. They have positive probability
	since $\nu_{t_2}$ has full support. By the law of total probability,
	\[
	\hat{\nu}_{t_3}(\cdot\mid\bar{\gs}=\xi)
	=\sum_{\rho:\,\bar{\rho}=\xi}
	\bP[\gs_{t_2}=\rho\mid\bar{\gs}_{t_2}=\xi]\,
	\bP[\hat{\gs}_{t_3}\in\cdot\mid\gs_{t_2}=\rho].
	\]
	Each law in the mixture is within $\gep/3$ of
	$\hat{\mu}_N(\cdot\mid\bar{\gs}=\xi)$. By convexity of total variation distance,
	\begin{equation}\label{step3:conditional}
	\bigl\|\hat{\nu}_{t_3}(\cdot\mid \bar{\gs}=\xi) - \hat{\mu}_N(\cdot\mid \bar{\gs}=\xi)\bigr\|_{\TV} \le \frac{\gep}{3}
	\end{equation}
	Moreover, $\bar{\nu}_{t_3}=\bar{\nu}_{t_2}$,
	so $\|\bar{\nu}_{t_3}-\bar{\mu}_N\|_{\TV}\le \gep/3$ by Proposition~\ref{prop:step2}.
	Conditioning on the skeleton $\bar{\gs}=\xi$ therefore gives
	\begin{align*}
	2\|\hat{\nu}_{t_3}-\hat{\mu}_N\|_{\TV}
	&= \sum_{\xi\in\bar{\cS}_N}\sum_{\hat{\gs}:\,\bar{\gs}=\xi}\bigl|\hat{\nu}_{t_3}(\hat{\gs})-\hat{\mu}_N(\hat{\gs})\bigr|\\
	&\le \sum_{\xi\in\bar{\cS}_N}\sum_{\hat{\gs}:\,\bar{\gs}=\xi}
	\Bigl[\bar{\nu}_{t_3}(\xi)\bigl|\hat{\nu}_{t_3}(\hat{\gs}\mid\bar{\gs}=\xi)-\hat{\mu}_N(\hat{\gs}\mid\bar{\gs}=\xi)\bigr|\\
	&\hspace{33mm} + \hat{\mu}_N(\hat{\gs}\mid\bar{\gs}=\xi)\bigl|\bar{\nu}_{t_3}(\xi)-\bar{\mu}_N(\xi)\bigr|\Bigr]\\
	&= 2\|\bar{\nu}_{t_3}-\bar{\mu}_N\|_{\TV}
	+ 2\sum_{\xi\in\bar{\cS}_N}\bar{\nu}_{t_3}(\xi)\bigl\|\hat{\nu}_{t_3}(\cdot\mid\bar{\gs}=\xi)-\hat{\mu}_N(\cdot\mid\bar{\gs}=\xi)\bigr\|_{\TV}
	\le \frac{4\gep}{3},
	\end{align*}
	Dividing by $2$ completes the proof.
\end{proof}

\subsection{Proof of Theorem~\ref{thm:main} and concluding remarks}
\begin{proof}[Proof of Theorem~\ref{thm:main}]
Fix $\gep\in(0,1)$. By Proposition~\ref{prop:upper} and Theorem~\ref{th:B}, for every
$\gd\in(0,1/300)$,
\[
\limsup_{N\to\infty}\frac{2\pi^2 t^N_{\Mix}(\gep)}{N^2\log N} \le 1+\gd.
\]
Letting $\gd\downarrow 0$ and combining the resulting bound with Proposition~\ref{prop:lower} yields
\begin{equation*}
1 \le \liminf_{N\to\infty}\frac{2\pi^2 t^N_{\Mix}(\gep)}{N^2\log N}
\le \limsup_{N\to\infty}\frac{2\pi^2 t^N_{\Mix}(\gep)}{N^2\log N} \le 1,
\end{equation*}
which proves~\eqref{eq:main}. Since the limit does not depend on $\gep\in(0,1)$, the
sequence of chains exhibits cutoff at time $(2\pi^2)^{-1}N^2\log N$.
\end{proof}

\begin{rem}[Cutoff window and spectral centering]\label{rem:window}
In view of~\eqref{gapinterchange} and~\eqref{eq:spectra}, Theorem~\ref{thm:main} can also be
stated as $t^N_{\Mix}(\gep) = (1+o(1))\log N/(2\gap_N)$, which implies
$t^N_{\Mix}(\gep)-t^N_{\Mix}(1-\gep) = o(N^2\log N)$ for every $\gep\in(0,1/2)$.
Because Assumption~\ref{assum:resist} specifies no quantitative rate of convergence for $\gd_N^{(0)}$,
the deterministic centering $(2\pi^2)^{-1}N^2\log N$ need not approximate $(2\gap_N)^{-1}\log N$ to
order $N^2$: for instance, when $r^{(N)}(x,x+1)=1+(\log N)^{-1/2}$ on every edge, the two scales
differ by $(1+o(1))(2\pi^2)^{-1}N^2\sqrt{\log N}$. This reflects the precision of the cutoff
location rather than the width of the cutoff window.
\end{rem}

\subsubsection*{Acknowledgements} The authors thank Hubert Lacoin for suggesting this problem, and thank Tertuliano Franco and Hubert Lacoin for enlightening discussions.
S.Y.\ is supported by CNPq Universal Grant 401314/2025-1, CNPq grant 408590/2024-6 and \mbox{FAPESP} grant
2025/27064-6. P.T.\ is supported by the  National Natural Science Foundation of China No.~12571151.
\subsubsection*{Statement on AI Use}
The authors used ChatGPT to assist with drafting and revising the exposition, improving readability, and clarifying and checking mathematical arguments. The authors take full responsibility for the accuracy and validity of the manuscript.

\end{document}